\documentclass[12pt]{article}
\usepackage{latexsym}
\usepackage{amssymb}
\begin{document}
\thispagestyle{empty}

\vskip 20pt
\begin{center}
{\bf NEW  LOWER BOUND FORMULAS FOR MULTICOLORED
RAMSEY NUMBERS}
\vskip 15pt
{\bf Aaron Robertson}\\
{\it Department of Mathematics,}
{\it Colgate University,
Hamilton, NY 13346}\\
{\tt aaron@math.colgate.edu}
\end{center}
\vskip 30pt
\begin{abstract}{\noindent
We give two lower bound formulas for
multicolored
Ramsey numbers.  These formulas improve
the bounds for several small multicolored Ramsey numbers.
}
\end{abstract} 

\vskip 20pt
{\bf 1. INTRODUCTION}

In this short note we give two new lower bound formulas
for the edgewise $r$-colored Ramsey numbers, $R(k_1,k_2,\dots,k_r)$.
Both formulas are derived via construction.

We will make use of the following notation.  Let $G$ be a graph, 
$V(G)$ the set of vertices of $G$, and $E(G)$ the set of edges of $G$.
An $r$-coloring, 
$\chi$, will be assumed to be an edgewise coloring, i.e. 
$\chi(G):E(G) \rightarrow \{1,2,\dots,r\}$.  If $u,v \in V(G)$, 
we take $\chi(u,v)$ to be the color of the edge
connecting $u$ and $v$ in $G$.
If we are considering the diagonal Ramsey numbers, 
i.e. $k_1=k_2=\dots,k_r=k$, we will use $R_r(k)$ to denote the
corresponding Ramsey number.  
It will also be helpful to make the following definition.

\noindent
{\bf Definition.}  A {\it Ramsey $r$-coloring for 
$R=R(k_1,k_2,\dots,k_r)$} is an $r$-coloring of the complete
graph on $V<R$ vertices which does not admit 
any monochromatic $K_{k_j}$ subgraph of color $j$ for $j = 1,2,\dots,r$.
For $V=R-1$ we call the coloring a {\it maximal Ramsey $r$-coloring}.

\section*{\normalsize 2. THE LOWER BOUNDS}

We start with a very trivial bound which nonetheless improves
upon some current best lower bounds.

\noindent
{\bf Theorem 1.}  {\it Let $r \geq 3$.  For any $k_i \geq 3$,
$i=1,2,\dots,r$, we have
$$R(k_1,k_2,\dots,k_r) > (k_1-1)(R(k_2,k_3,\dots, k_r)-1).$$}

\noindent
{\it Proof.}  Let $\phi(G)$ be a maximal Ramsey $(r-1)$-coloring for
$R(k_2,k_3,\dots,k_r)$ with colors $2,3,\dots,r$.  Let $k_1 \geq 3$.
Let $G_i=G$, for $i=1,2,\dots,k_1-1$.  Let $v_i \in G_i$, 
$v_j \in G_j$ and define $\chi(H)$ as follows:
$$
\chi(v_i,v_j) = \left\{
\begin{array}{ll}
\phi(v_i,v_j)&\mbox{if} \,\, i=j\\
1&\mbox{if} \,\, i \neq j.\\
\end{array}
\right.
$$

We now show that $\chi(H)$ is a Ramsey $r$-coloring for
$R(k_1,k_2,\dots,k_r)$.  For $j \in \{2,3,\dots,r\}$, 
$\chi(H)$ does not admit any monochromatic $K_{k_j}$ of color $j$ by
the definition of $\phi$.  Hence, we need only consider
color $1$.  Since $\phi(G_i)$, $1 \leq i \leq k_1-1$, 
is void of color $1$, any monochromatic $K_{k_1}$ of color $1$ may only
have one vertex in $G_i$ for each $i \in \{1,2,\dots,k_1-1\}$.
By the pigeonhole principle, however, there
exists $I \in \{1,2,\dots,k_1-1\}$ such that $G_I$
contains two vertices of $K_{k_j}$, a contradiction.
\hfill{$\Box$}

\noindent
{\it Examples.}  Theorem 1 implies that
$R_5(4) \geq 1372, R_5(5) \geq 7329, R_4(6) \geq 5346,$ and $R_4(7) \geq
19261$, all of which beat the current best known bounds given
in [Rad].

We now look at an off-diagonal bound.

\noindent
{\bf Theorem 2.}  {\it  Let $r \geq 3$.  For any $3 \leq k_1<k_2$, and
$k_j
\geq 3$,
$j =3,4,\dots,r$ we have
$$R(k_1,k_2,\dots,k_r) > (k_1+1)(R(k_2-k_1+1,k_3,\dots,
k_r)-1).$$}

Before giving the proof of this theorem, we have need of
the following definition.

\noindent
{\bf Defintion.}  We say that the $n \times n$ symmetric matrix
$T=T(x_0,x_1,\dots,x_r)$ is a {\it Ramsey incidence matrix} if
the $r$-coloring defined by 
$\chi:E(K_n) \rightarrow \{x_1,x_2,\dots,x_r\}$, 
$\chi(i,j)=(i,j)$, is
a Ramsey $r$-coloring.
Furthermore, the color $x_0$ appears only on the diagonal
of $T$ (which we will
denote $diag(T)$).  Note that $T(x_0,x_2,x_1,x_3,\dots,x_r)$ defines
the same graph as $T(x_0,x_1,x_2,x_3,\dots,x_r)$ with
colors $x_1$ and $x_2$ interchanged.

\noindent 
{\it Proof of Theorem 2.}  We will construct an $r$-colored complete graph
on
$(k_1+1)(R(k_2-k_1+1,k_3,\dots,k_r)-1)$ vertices which avoids 
monochromatic subgraphs $K_{k_i}$ of color $i$, $i=1,2,\dots,r$, by
means of Ramsey incidence matrices.  We start the proof with
$R(t,k,l)$ and then generalize to an arbitrary number of colors.  

Consider a maximal Ramsey $2$-coloring for
$R(k,l-t+1)$.  Let $T=T(x_0,x_1,x_2)$
denote the associated Ramsey incidence matrix.  Define 
$A:=T(\bullet,2,3)$, $B:=T(3,2,1)$, and $C:=T(1,2,3)$ and
consider the symmetric $(t+1) \times (t+1)$ block matrix, $H$, below.

$$
\begin{array}{ccccccccc}
&A\\
&B&A\\
&C&C&A\\
H=&C&C&B&A\\
&C&C&B&B&A\\
&\vdots&\vdots&\vdots&\vdots&\ddots&\ddots\\
&C&C&B&B&\dots&B&A
\end{array}
$$

We will show that $H$ contains no monochromatic $K_t$ of color $1$, no
monochromatic $K_k$ of color $2$, and no monochromatic 
$K_l$ of color $3$,
for
$l \geq t+1$, to show that $R(t,k,l) > (t+1)(R(k,l-t+1)-1)$.
To this end, we
first look at the stucture of the edges of $K_{s}$ in $H$.
Without loss of generality we may assume that the entries in $H$
representing the edges of $K_{s}$ have the following structure, 
where $j_1 < i_1$.

$$
\begin{array}{llllll}
\bullet(i_1,j_1)\\
\bullet(i_2,j_1)&\bullet (i_2,i_1)\\
\bullet(i_3,j_1)&\bullet (i_3,i_1)&\bullet (i_3,i_2)\\
&\hskip -40pt\vdots&\hskip -39pt\vdots \hskip 55pt \vdots&\ddots\\
\bullet(i_{s-1},j_1)&\bullet(i_{s-1},i_1)&\bullet
(i_s,i_2)&\dots&\bullet (i_{s-1},i_{s-2})\\
\end{array}
$$

We will refer to
two different types of rows below:  entry rows and block rows.
An entry row is a set $\{(i_q,j_r): 1 \leq r \leq q\}$ were
$q$ is a fixed integer between $1$ and $s-1$.  A block row consists of one
of the rows of $H$, for example the third block
row is $C C A$.  We will also use the
term {\it relative position of $K_s$} several times.  To determine the
relative position of $K_s$, take all of the corresponding
coordinates of $K_s$ in H and reduce them
modulo $(R(t,k,l)-1)$.  This reduction gives us entries only in the
$A(1,1)$ block in $H$.  (When confusion may arise, we will use the full
notation
$A(i,j)$ to clarify which $A$ (or $B$ or $C$) block is being considered).

We will now show that the graph defined by the Ramsey incidence
matrix $H$ avoids the desired monochromatic subgraphs.

\noindent
{\bf No monochromatic $\mathbf{K_t}$ of color $\mathbf{1}$}.  $K_t$ cannot
have two entry rows in any block row containing a $C$ since $1 \in
diag(C)$ and $1 \not \in A$.
Further, if two block rows both containing a $C$ have entry rows in them,
then since $1 \in diag(C)$ we must have $1 \in diag(B)$, a 
contradiction.  Hence, $B(2,1)$ must have at least two entry
rows.  This implies that $1 \in A(2,2)$, a contradiction.  Thus,
we cannot have a monochromatic $K_t$ of color 1.

\noindent
{\bf No monochromatic $\mathbf{K_k}$ of color $\mathbf{2}$}.  If a
monochromatic
$K_k$ of color 2 exists in $H$, then by taking the relative position, we
would have a monochromatic $K_k$ of color 2 in $A(1,1)$, contradicting the
definition of $A$.

\noindent
{\bf No monochromatic $\mathbf{K_l}$ of color $\mathbf{3}$}. 
Assume, for a contradiction, that a monochromatic
$K_l$ of color 3 exists.  If
there are no entries in any $B$, then taking the relative position
of all entries will
imply that $A(1,1)$ contains a monochromatic $K_l$ of color 3, a
contradiction. Hence, we must have at least one entry in some $B$.
However, each $B$ may contain at most one entry since $3 \in diag(B)$.
This implies that we can have at most one entry row in each block
row $4$ through $t+1$, and at most one entry column in the
first block column.  We now delete the first entry column,
and the bottom $t-2$ entry rows.  This deletion procedure
assures us that none of the remaining entries lie in any $B$.  
Hence, we are left with $l-t$ entry rows, 
which form a $K_{l-t+1}$.  By taking the relative position of these
remaining entries,  we have a monochromatic $K_{l-t+1}$ of color 3
in $A(1,1)$, a contradiction.

To generalize the above argument to an arbitrary number of
colors we change the definitions of $A$, $B$, and $C$;
$A:=T(0,2,3,4,5,\dots,r)$,
$B:=T(3,2,1,4,5,\dots,r)$, $C:=T(1,2,3,4,5,\dots,r)$.  To see that
there is no monochromatic $K_{k_j}$ of color $j$ for $j=4,5,\dots,r$ see
the argument for no monochromatic $K_k$ of color 2 above.
\hfill{$\Box$}

\noindent
{\it Example.}  Theorem 2 implies that
$R(3,3,3,11) \geq 437$, beating the previous best
lower bound of $433$.

\section*{\normalsize REFERENCES}

\noindent
[Rad] S. Radziszowski, Small Ramsey numbers, {\it El.
J. Comb.}, DS1 (revision \#7, 2000), 36pp.

\noindent
[Ram] F. Ramsey, On a problem of formal logic,
{\it Proc. London Math. Soc.} {\bf 30} (1930), 264-286.

\noindent
[Rob] A. Robertson, {\it Ph.D. thesis}, Temple University, 1999.

\end{document}